\documentclass{amsart}

\usepackage{fullpage}
\usepackage{amsmath}
\usepackage{amsfonts}
\usepackage{amssymb}

\input xy
\xyoption{matrix}\xyoption{arrow}\xyoption{curve}  \xyoption{frame}

\begin{document}

%%%%% Functors %%%%%

\newcommand{\Hom}{\mathrm{Hom}}
\newcommand{\RHom}{\mathrm{RHom}^*}
\newcommand{\HOM}{\mathrm{HOM}}
\newcommand{\stHom}{\underline{\mathrm{Hom}}}
\newcommand{\Ext}{\mathrm{Ext}}
\newcommand{\Tor}{\mathrm{Tor}}
\newcommand{\HH}{\mathrm{HH}}
\newcommand{\Endo}{\mathrm{End}}
\newcommand{\ENDO}{\mathrm{END}}
\newcommand{\stEnd}{\mathrm{\underline{End}}}
\newcommand{\Tr}{\mathrm{Tr}}

%%%%% Functions/Operators (no space after) %%%

\newcommand{\coker}{\mathrm{coker}}
\newcommand{\aut}{\mathrm{Aut}}
\newcommand{\op}{\mathrm{op}}
\newcommand{\ob}{\mathrm{ob}}

\newcommand{\add}{\mathrm{add}}
\newcommand{\ADD}{\mathrm{ADD}}
\newcommand{\ind}{\mathrm{ind}}
\newcommand{\rad}{\mathrm{rad}}
\newcommand{\soc}{\mathrm{soc}}
\newcommand{\ann}{\mathrm{ann}}
\newcommand{\im}{\mathrm{im}}
\newcommand{\chr}{\mathrm{char}}
\newcommand{\pdim}{\mathrm{p.dim}}

%%%%% Categories %%%%%%%%%

\newcommand{\rmod}{\mbox{mod-}}
\newcommand{\Rmod}{\mbox{Mod-}}
\newcommand{\lmod}{\mbox{-mod}}
\newcommand{\lMod}{\mbox{-Mod}}
\newcommand{\stmod}{\mbox{\underline{mod}-}}
\newcommand{\stlmod}{\mbox{-\underline{mod}}}

\newcommand{\gr}{\mbox{gr-}}
\newcommand{\stgr}{\mbox{\underline{gr}-}}

\newcommand{\gmod}[1]{\mbox{mod}_{#1}\mbox{-}}
\newcommand{\gMod}[1]{\mbox{Mod}_{#1}\mbox{-}}
\newcommand{\Bimod}[1]{\mathrm{Bimod}_{#1}\mbox{-}}

\newcommand{\proj}{\mbox{proj-}}
\newcommand{\lproj}{\mbox{-proj}}
\newcommand{\Proj}{\mbox{Proj-}}
\newcommand{\inj}{\mbox{inj-}}
\newcommand{\coh}{\mbox{coh-}}

%%%%% Style %%%%%%%%%%

\newcommand{\und}[1]{\underline{#1}}
\newcommand{\gen}[1]{\langle #1 \rangle}
\newcommand{\floor}[1]{\lfloor #1 \rfloor}
\newcommand{\ceil}[1]{\lceil #1 \rceil}
\newcommand{\bnc}[2]{ \left( \scriptsize{\begin{array}{c} #1 \\ #2 \end{array} } \right)}
\newcommand{\bimo}[1]{{}_{#1}#1_{#1}}
\newcommand{\ses}[5]{\ensuremath{0 \rightarrow #1 \stackrel{#4}{\longrightarrow} 
#2 \stackrel{#5}{\longrightarrow} #3 \rightarrow 0}}
\newcommand{\bfa}{\mathbf{a}}
\newcommand{\bfb}{\mathbf{b}}
\newcommand{\bfc}{\mathbf{c}}
\newcommand{\bfd}{\mathbf{d}}
\newcommand{\A}{\mathcal{A}}
\newcommand{\B}{\mathcal{B}}
\newcommand{\C}{\mathcal{C}}
\newcommand{\D}{\mathcal{D}}
\newcommand{\E}{\mathcal{E}}
\newcommand{\F}{\mathcal{F}}
\newcommand{\N}{\mathcal{N}}
\newcommand{\T}{\mathcal{T}}
\newcommand{\s}{\mathcal{S}}
\newcommand{\X}{\mathcal{X}}
\newcommand{\Y}{\mathcal{Y}}
\newcommand{\ppX}{{}^{\perp}\mathcal{X}^{\perp}}
\newcommand{\tC}{\tilde{\mathcal{C}}}
\newcommand{\tK}{\tilde{K}(A)}
\newcommand{\ul}[1]{\underline{#1}}

%%%%% Theorems %%%%%%%%%%

\newtheorem{therm}{Theorem}[section]
\newtheorem{defin}[therm]{Definition}
\newtheorem{propos}[therm]{Proposition}
\newtheorem{lemma}[therm]{Lemma}
\newtheorem{coro}[therm]{Corollary}

\title{Torsion pairs and simple-minded systems in triangulated categories}
\author{Alex Dugas}
%\affil{University of the Pacific}

\address{Department of Mathematics, University of the Pacific, 3601 Pacific Ave, Stockton CA 95211, USA.}
\email{adugas@pacific.edu} 
\thanks{\noindent Corresponding Author: Alex Dugas, University of the Pacific, 3601 Pacific Ave, Stockton CA 95211, USA\\ Email:  {\tt adugas@pacific.edu}}

%\subjclass[2010]{16G10, 16E35, 18E30}
%\keywords{derived equivalence, stable equivalence, quiver mutation, weakly symmetric algebra}

\begin{abstract} 
Let $\mathcal{T}$ be a Hom-finite triangulated Krull-Schmidt category over a field $k$.  Inspired by a definition of Koenig and Liu \cite{KoeLiu}, we say that a family $\mathcal{S} \subseteq \mathcal{T}$ of pairwise orthogonal bricks is a simple-minded system if its closure under extensions is all of $\mathcal{T}$.  We construct torsion pairs in $\T$ associated to any subset $\mathcal{X}$ of a simple-minded system $\s$, and use these to define left and right mutations of $\s$ relative to $\X$.  When $\T$ has a Serre functor $\nu$ and $\s$ and $\X$ are invariant under $\nu \circ [1]$, we show that these mutations are again simple-minded systems.  We are particularly interested in the case where $\T = \stmod \Lambda$ for a self-injective algebra $\Lambda$.  In this case, our mutation procedure parallels that introduced by Koenig and Yang for simple-minded collections in $D^b(\rmod \Lambda)$ \cite{KoYa}.  It follows that the mutation of the set of simple $\Lambda$-modules relative to $\X$ yields the images of the simple $\Gamma$-modules under a stable equivalence $\stmod \Gamma \rightarrow \stmod \Lambda$, where $\Gamma$ is the tilting mutation of $\Lambda$ relative to $\X$.  
\end{abstract}

\maketitle

 \section{Introduction} 
\setcounter{equation}{0}

A well-known conjecture of Auslander and Reiten states that two finite-dimensional algebras with equivalent stable module categories have the same number of nonprojective simple modules (up to isomorphism) \cite{ARS}, and work of Mart\'inez-Villa has reduced this problem to the case where the algebras are self-injective \cite{IPUSE}.  One possible approach to this conjecture involves studying the sets of objects in the stable category $\stmod \Lambda$ of a finite-dimensional self-injective algebra $\Lambda$, which correspond to the set of simple $\Gamma$-modules under some equivalence $\stmod \Gamma \stackrel{\approx}{\longrightarrow} \stmod \Lambda$ for another finite-dimensional self-injective algebra $\Gamma$.  If one can provide a characterization of such sets that is intrinsic to $\stmod \Lambda$, then perhaps it is possible to prove that all such sets have the same cardinality.  In this direction, Pogorza\l y has taken a few obvious properties of such sets of objects in $\stmod \Lambda$ as the basis for his definition of {\it maximal systems of stable orthogonal bricks} in \cite{Pog}.  More recently Koenig and Liu have refined this notion in their introduction of {\it simple-minded systems} \cite{KoeLiu} (although, to the best of our knowledge, it is not yet known whether these two notions are truly distinct).  Briefly, a set of pairwise orthogonal indecomposable objects in $\stmod \Lambda$ with endomorphism rings equal to division rings forms a simple-minded system if it generates $\stmod \Lambda$ by extensions.  This latter condition is naturally difficult to check, and thus non-trivial examples of simple-minded systems are hard to find, especially in the absence of a non-trivial stable equivalence.  The primary goal of this article is to define a mutation procedure for simple-minded systems that allows one to easily construct many non-trivial examples starting out from the set of simple modules itself.

Our mutation procedure is a generalization of that introduced in \cite{MutSym}, where it was defined only for maximal systems of stable orthogonal bricks and in a more restricted setting.  Like our previous definition, it is inspired by Okuyama's description of the images of the simple modules under a specific stable equivalence from a {\it tilting mutation} of $\Lambda$.  More precisely, for a symmetric algebra $\Lambda$ and the endomorphism ring $\Gamma$ of an Okuyama tilting complex $T$ associated to a subset $U$ of the simple $\Lambda$-modules, Okuyama described the images of the simple $\Gamma$-modules under the stable equivalence $\ul{F}$ induced by the derived equivalence $D^b(\rmod \Gamma) \rightarrow D^b(\rmod \Lambda)$ given by $T$.   Writing $S_i$ for the simple $\Lambda$-modules and $S'_i$ for the simple $\Gamma$-modules ($1 \leq i \leq n$), Okuyama's Lemma (Lemma 2.1' in \cite{Oku}) states that 
\begin{itemize} \item For each $S_i \in U$,  $\ul{F}(S'_i) \cong S_i$; and 
\item For each $S_j \notin U$, $\ul{F}(S'_j)$ is isomorphic to the largest submodule $X$ of $\Omega(S_j)$ with $\Hom_{\Lambda}(X,S_i) = 0$ for all $S_i \in U$.  
\end{itemize}

The images of the simple $\Gamma$-modules under this stable equivalence naturally form a simple-minded system in $\stmod \Lambda$, which coincides with our definition of the mutation $\mu^+_U(\s)$ of the set $\s$ of simple $\Lambda$-modules relative to $U$.  By describing this transformation taking $S_i$ to $\ul{F}(S'_i)$ in terms of the triangulated structure of $\stmod \Lambda$, we will obtain our definition of mutation for arbitrary simple-minded systems in a triangulated category (see Definition 4.1).  In particular, it is possible to apply such mutations in succession, and thereby track, inside $\stmod \Lambda$, the various derived equivalences to $D^b(\rmod \Lambda)$ that can be obtained by successively tilting via Okuyama complexes.

This paper is organized as follows.  We begin working in an abstract Hom-finite triangulated Krull-Schmidt category $\T$ over a fixed field $k$.  In Section 2, we develop some basic properties of filtration subcategories generated by a maximal system of orthogonal bricks in $\T$.  As our definition of mutation of simple-minded systems makes use of torsion pairs in $\T$, we study these in Section 3.  We show, in particular, that one obtains natural examples corresponding to any subset of a simple minded system in $\T$.  Using these torsion pairs, we then define left and right mutations of a simple-minded system $\s$ at a subset $\X$ in Section 4.  Here, we assume that $\T$ has a Serre functor $\nu$ and that $\s$ and $\X$ are stable under $\nu \circ [1]$.  Our main result states that these mutations are again simple-minded systems.  The bulk of the proof involves showing that the mutation of a simple-minded system again generates $\T$ by extensions.  It turns out that in a stable module category $\stmod \Lambda$ our mutations of simple-minded systems closely parallel the mutations of {\it simple-minded collections} in the derived $D^b(\rmod \Lambda)$ as studied by Koenig and Yang in \cite{KoYa}.  Thus, in Section 5 we describe how new results of Koenig and Yang connect mutation of simple minded systems in $\stmod \Lambda$ to tilting mutations of $\Lambda$ and to {\it silting mutation} in $K^b(\proj \Lambda)$ as introduced by Aihara and Iyama \cite{AihIya}.  In particular, we show how Okuyama's Lemma follows from this work of Koenig and Yang.  Finally, we consider several examples in the last section.

Throughout this article, $k$ will denote a field and all categories and functors are assumed to be $k$-linear.  We will use $\Lambda$ and $\Gamma$ to denote finite-dimensional, basic, self-injective $k$-algebras.  We let $\rmod \Lambda$ denote the category of finitely presented right $\Lambda$-modules, and $\stmod \Lambda$ the associated stable category obtained by factoring out the ideal of maps that factor through a projective.  We also consider the homotopy category $K^b(\proj \Lambda)$ of bounded complexes of finitely-generated projective $\Lambda$-modules and the derived category $D^b(\rmod \Lambda)$.  Unless otherwise noted, all modules considered are right modules, and we write morphisms on the left, composing them from right to left.  The same convention is used for morphisms in an arbitrary category.    For a category $\C$, we will also write $\C(X,Y)$ for the set of morphisms between two objects $X$ and $Y$ of $\C$.

\section{Filtration subcategories}
\setcounter{equation}{0}

Let $k$ be a field and let $\T$ be a Hom-finite Krull-Schmidt triangulated $k$-category with suspension denoted by $[1]$.  We write $\C \subseteq \T$ for a subset of $\mbox{ob}(\T)$, which we identify with the corresponding full subcategory of $\T$.  For $\C, \D \subseteq \T$, we set $\gen{\C} = \add(\C)$ and $$\C *\D = \{X \in \T\ |\ \mbox{there\ is\ a\ triangle\ } C \rightarrow X \rightarrow D \rightarrow\ \mbox{with}\ C \in \C, D \in \D\}.$$
The next lemma shows that formation of these extension subcategories is associative.

\begin{lemma}[Cf. Lemma 1.3.10 in \cite{BBD}] Let $\C, \D, \E \subseteq \T$.  Then $(\C * \D) * \E = \C * (\D * \E)$.
\end{lemma}

\noindent
{\it Proof.}  Let $X \in (\C * \D) * \E$ so that we have triangles $X' \longrightarrow X \longrightarrow Y \rightarrow$ and $Z \longrightarrow X' \longrightarrow W \rightarrow$ with $Z \in \C, W \in \D$ and $Y \in \E$.  By the octahedral axiom we have a commutative diagram where the rows and columns are distinguished triangles in $\T$.
$$\xymatrix{Z \ar[d] \ar@{=}[r] & Z \ar[d] \\ X' \ar[d] \ar[r] & X \ar[d] \ar[r] & Y \ar[r] \ar@{=}[d] & \\ W \ar[d] \ar[r] & U \ar[d] \ar[r] & Y \ar[r]  & \\ & }$$
The bottom row shows that $U \in \D * \E$, and now the second column shows that $X \in \C * (\D *\E)$.  The reverse inclusion is proved similarly.  $\Box$ \\

We say that $\C$ is {\it extension-closed} if $\C*\C \subseteq \C$.  For $\C \subseteq \T$, we can define various filtration subcategories inductively as follows:

\begin{eqnarray} (\C)_n & = & \left\{ \begin{array}{cl} \{0\} & \mbox{if}\ n=0 \\ (\C)_{n-1} * (\C \cup \{0\}) & \mbox{if}\ n \geq 1 \end{array} \right. \\
{}_n(\C) & = & \left\{ \begin{array}{cl} \{0\} & \mbox{if}\ n=0 \\  (\C \cup \{0\}) * {}_{n-1}(\C) & \mbox{if}\ n \geq 1 \end{array} \right.
\end{eqnarray}

Observe that these filtration subcategories form ascending chains: $\{0\} = (\C)_0 \subseteq (\C)_1 \subseteq (\C)_2 \subseteq \cdots$ and similarly for the ${}_n(\C)$.  Furthermore, observe that each $(\C)_n$ and each ${}_n(\C)$ is a {\it strict} subcategory, meaning that it is closed under isomorphisms.  However, in general, these subcategories need not be closed under direct summands.  We begin by showing that these two chains of subcategories are identical.  

\begin{lemma} For any $\C \subseteq \T$ we have $(\C)_n = {}_n(\C)$ for all $n \geq 0$.
\end{lemma}

\noindent
{\it Proof.}  Clearly the equality holds for $n=0, 1$.   Using Lemma 2.1 it is now easy to see that it holds for all $n \geq 0$ by induction.  $\Box$ \\

We now define the {\it filtration subcategory generated by} $\C$ to be $\F(\C) = \cup_{n \geq 0} (\C)_n$.   It follows from the next lemma that $\F(\C)$ is the smallest extension-closed subcategory of $\T$ containing $\C$.

\begin{lemma} For any $\C \subseteq \T$, we have $(\C)_n * (\C)_m = (\C)_{n+m}$ for all $m, n \geq 0$.  In particular, the filtration subcategory $\F(\C)$ is closed under extensions.
\end{lemma}

\noindent
{\it Proof.}  We use induction on $m$.  By definition of $(\C)_n$, the equality holds for all $n \geq 0$ when $m=0, 1$.  Now suppose that it holds for all $n \geq 0$ and all $k < m$, for some $m \geq 2$.  By Lemma 2.1, we have $(\C)_n*(\C)_m = ((\C)_n*(\C)_{m-1})*(\C \cup \{0\}) = (\C)_{n+m-1}*(\C \cup \{0\}) = (\C)_{n+m}$. $\Box$ \\

We are mainly interested in filtration subcategories that are generated by a set of objects that behave like simple objects in $\T$.  The following definitions, inspired by \cite{Pog}, are slightly more general than those used in \cite{MutSym}.  

\begin{defin} A set $\s = \{S_i\}_{i \in I}$ of objects of $\T$ is set of {\bf (pairwise) orthogonal bricks} in $\T$ if $$ (1)\ \ \ \T(S_i,S_j) \cong \left\{ \begin{array}{rl} 0, & \mbox{if\ }i \neq j \\ \mbox{a\ division\ ring}\ k_i, & \mbox{if} \ i=j. \end{array} \right.$$
  We say that $\s$ is a {\bf maximal system of orthogonal bricks} in $\T$ if additionally  \begin{eqnarray*} (2)\ & & \forall\  X \in \T\ \exists\ i \in I\ \mbox{such\ that}\ \T(X,S_i) \neq 0; and \\ 
(2')\ & & \forall\ Y \in \T\ \exists\ i \in I\ \mbox{such\ that}\ \T(S_i, Y) \neq 0.\end{eqnarray*}   
\end{defin}

We will also consider the following stronger definition, inspired by Koenig and Liu \cite{KoeLiu}.

\begin{defin} A set $\s$ of orthogonal bricks in $\T$ is a {\bf simple-minded system} if $\F(\s) = \T$.  In this case, any object $X \in \T$ has an {\bf $\s$-length} $l_\s(X) := \min \{n \geq 0\ |\ X \in (\s)_n\}$.
\end{defin}

\noindent
{\it Remarks.}  (1) As in \cite{KoeLiu}, it is not hard to see that a simple-minded system is a maximal system of orthogonal bricks.  We do not know if the converse holds in general. 

(2) If $\T$ has a Serre functor $\nu$ inducing a Serre duality $\T(X,Y) \cong D\T(Y, \nu X)$ for all $X, Y \in \T$, where $D = \Hom_k(-,k)$ is the standard $k$-duality, then (2) and (2') in the definition of maximal system of orthogonal bricks are easily seen to be equivalent.  Moreover, in the presence of a Serre duality, it is natural to require a maximal system of orthogonal bricks, or  a simple-minded system, $\s$ to satisfy the additional condition $$(3) \hspace{2cm}  \nu(\s[1]) = \s. \hspace{2cm}$$  This is due to the fact that if $\T = \stmod A$ for a self-injective algebra $A$, then $\nu = \N \Omega$ so that $\nu(X[1]) = \N(X)$, where $\N$ denotes the Nakayama functor.  In particular, since $\N$ permutes the simple $A$-modules, and commutes with any stable equivalence, any simple-minded system in $\stmod B$ which arises as the image of the simple $A$-modules under an equivalence $\stmod A \rightarrow \stmod B$ will be stable under $\N$.\\

\begin{lemma}  Let $\s \subseteq \T$ be a set of orthogonal bricks, and suppose that there is a distinguished triangle $Y \longrightarrow X \stackrel{f}{\longrightarrow} S_i \rightarrow$ with $X \in (\s)_n$ for some $n \geq 1$, $S_i \in \s$ and $f \neq 0$.  Then $Y \in (\s)_{n-1}$.
\end{lemma}

\noindent
{\it Proof.}  We proceed by induction on $n \geq 1$.  If $n = 1$, then $f \neq 0$ implies that $f$ is an isomorphism.  Hence $Y = 0 \in (\s)_0$.  Now assume that the statement is true for $n-1$, where $n \geq 2$.  Since $X \in (\s)_n$, we have a triangle $X' \longrightarrow X \stackrel{g}{\longrightarrow} S_j \rightarrow$ with $S_j \in \s$ and $X' \in (\s)_{n-1}$.  Note that we may assume that $S_j \neq 0$, as otherwise we would have $X \cong X' \in (\s)_{n-1}$, and we would be done by the induction hypothesis.  By the octahedral axiom we have a commutative diagram where the rows and columns are distinguished triangles in $\T$.
$$\xymatrix{& X' \ar[d]^i \ar@{=}[r] & X' \ar[d]^{fi} \\ Y \ar@{=}[d] \ar[r] & X \ar[d]^{g} \ar[r]^{f} & S_i \ar[r] \ar[d] & \\ Y  \ar[r] & S_j \ar[d] \ar[r] & Z \ar[r]  \ar[d] & \\ & & }$$
By the rotation axiom, we thus obtain triangles $Z[-1] \longrightarrow X' \stackrel{fi}{\longrightarrow} S_i \rightarrow$ and $Z[-1] \longrightarrow Y \longrightarrow S_j \rightarrow$.  As long as $fi \neq 0$, we can use the induction hypothesis to conclude that $Z[-1] \in (\s)_{n-2}$ and then that $Y \in (\s)_{n-1}$.  Thus assume that $fi = 0$.  Thus $f$ factors as $f = hg$ for some nonzero $h : S_j \rightarrow S_i$.  Hence $h$ must be an isomorphism and it follows that $Y \cong X' \in (\s)_{n-1}$.  $\Box$ \\

\begin{lemma} If $\s \subseteq \T$ is a set of orthogonal bricks, then $(\s)_n$ is closed under direct summands for each $n \geq 1$.  In particular, $\F(\s)$ is closed under direct summands.
\end{lemma}

\noindent
{\it Proof.}  We proceed by induction on $n$.  For $n=1$ the statement is clear since each $S_i \in \s$ is indecomposable.  Thus let $n \geq 2$ and assume that $(\s)_{n-1}$ is closed under direct summands.  Suppose that $X \oplus Y \in (\s)_n$, and that we have a triangle $X' \stackrel{\bnc{a}{b}}{\longrightarrow} X \oplus Y \stackrel{(c\ d)}{\longrightarrow} S_i \rightarrow$ with $S_i \in \s$ and $X' \in (\s)_{n-1}$.  If $c = 0$ then the inclusion $X \rightarrow X \oplus Y$ factors through $\bnc{a}{b}$, which implies that $1_X$ factors through $a$.  Thus $X$ is isomorphic to a direct summand of $X'$ and so $X \in (\s)_{n-1} \subseteq (\s)_{n}$ by the induction hypothesis.  Now assume $c \neq 0$.  We thus have triangles $X'' \longrightarrow X \stackrel{c}{\longrightarrow} S_i \rightarrow$ and $X'' \oplus Y \longrightarrow X \oplus Y \stackrel{(c\ 0)}{\longrightarrow} S_i$, with $X'' \oplus Y \in (\s)_{n-1}$ by Lemma 2.6.  By the induction hypothesis, we obtain $X'' \in (\s)_{n-1}$ and hence $X \in (\s)_n$.  $\Box$ \\

\section{Torsion Pairs}
\setcounter{equation}{0}

In this section we study torsion pairs in triangulated categories as introduced by Iyama and Yoshino \cite{IY}.  This definition is slighty weaker than that used by Beligiannis and Reiten \cite{BelRei}, making it better suited for stable module categories than for derived categories.  For $\C \subseteq \T$ we set 
\begin{eqnarray*} \C^{\perp} & = & \{ X \in \T\ |\ \T(C, X) = 0\ \mbox{for\ all}\ C \in \C\} \\ {}^{\perp}\C & = & \{X \in \T\ |\ \T(X,C)=0\ \mbox{for\ all\ } C \in \C\}
\end{eqnarray*}
It is immediate that both $\C^{\perp}$ and ${}^{\perp}\C$ are strict, extension closed subcategories of $\T$ which are also closed under direct summands.  We shall also write ${}^{\perp}\C^{\perp}$ for the intersection ${}^{\perp}\C \cap \C^{\perp}$.

\begin{defin}  A pair $(\C,\D)$ of  full, additive subcategories of $ \T$, which are closed under direct summands, form a {\bf torsion pair}  if the following hold:
\begin{itemize}
\item[(TP1)] $\T(\C,\D) = 0$;
\item[(TP2)] $\T = \C*\D$.  That is, for each $X \in \T$, there exists a distinguished triangle $C_X \stackrel{f_X}{\longrightarrow} X \stackrel{g_X}{\longrightarrow} D_X \rightarrow$ with $C_X \in \C$ and $D_X \in \D$.
\end{itemize}
\end{defin}

In light of (TP2), we easily see that (TP1) can be strengthened to $$\mbox{(TP1')}\ \ \C^{\perp} = \D\ \ \mbox{and} \ \ \ \C = {}^{\perp}\D.$$  In particular, we see that the subcategories making up a torsion pair are necessarily extension-closed.  In any triangle as in (TP2), notice that (TP1) implies that $f_X$ is a right $\C$-approximation and $g_X$ is a left $\D$-approximation.  In particular, if $(\C,\D)$ is a torsion pair then $\C$ is contravariantly finite in $\T$ and $\D$ is covariantly finite in $\T$.  It is natural to wonder whether this triangle is unique up to isomorphism.  A priori, there is no reason it must be.  However, since $\T$ is assumed to be a Krull-Schmidt category, it is possible to choose a right minimal version of $f_X$ and then the resulting triangle is unique up to isomorphism (see 4.2-4.3 in \cite{GSE} for information about right and left minimal morphisms in a Krull-Schmidt category).

\begin{lemma} Let $(\C, \D)$ be a torsion pair in $\T$, and consider a triangle $C \stackrel{f}{\longrightarrow} X \stackrel{g}{\longrightarrow} D\stackrel{h}{\rightarrow}$.  The following are equivalent
\begin{enumerate} \item $C \in \C$ and $f$ is a minimal right $\C$-approximation;  
\item $D \in \D$ and $g$ is a minimal left $\D$-approximation.
\end{enumerate}
\end{lemma}

\noindent
{\it Proof.}  $(1) \Rightarrow (2)$: Since $(\C,\D)$ is a torsion pair, we have a triangle $C_X \stackrel{f_X}{\longrightarrow} X \stackrel{g_X}{\longrightarrow} D_X \rightarrow$ with $C_X \in \C$ and $D_X \in \D$.  As $f_X$ and $f$ are both right $\C$-approximations, we must have $C_X \cong C \oplus C'$ with $f_X|_{C'} = 0$.  It follows that $D_X \cong D \oplus C'[1]$, whence $D \in \D$.  Furthermore, if $g$ is not left minimal, $D$ contains a direct summand $D'$ on which $h$ restricts to an isomorphism.  It follows that $C$ has a direct summand, isomorphic to $D'[-1]$ on which $f$ restricts to $0$, but this contradicts the minimality of $f$.  The converse is proved analogously.  $\Box$ \\

We call a triangle  of the form $C \stackrel{f}{\longrightarrow} X \stackrel{g}{\longrightarrow} D \rightarrow$ with $f$ a minimal right $\C$-approximation a {\it minimal $(\C,\D)$-triangle} for $X$.  By choosing a right minimal version of $f_X$ in the triangle in (TP2), we know that any $X$ admits a minimal $(\C,\D)$-triangle.  Moreover, since a minimal right $\C$-approximation is unique up to isomorphism, so is the resulting triangle.  Furthermore, notice that a $(\C,\D)$-triangle is necessarily minimal if  $C$ (resp.  $D$) is indecomposable and $f$ (resp. $g$) is nonzero.

Fixing a torsion pair $(\C,\D)$ in $\T$, we can now define operators $\bfa$ and $\bfb$ on $\ob(T)$ via a minimal $(\C,\D)$-triangle $C_X \stackrel{f_X}{\longrightarrow} X \stackrel{g_X}{\longrightarrow} D_X \rightarrow$ for $X$ by $\bfa(X) = C_X$ and $\bfb(X) = D_X$.  These are not functors in general, but in certain cases they become functorial if we restrict to an appropriate subcategory (see the remarks following Lemma 4.7).  Notice, additionally, that they are functorial if one assumes $\C[1] \subseteq \C$ and $\D[-1] \subseteq \D$ as in the definition of torsion pair given by Beligiannis and Reiten \cite{BelRei}.  Moreover, under these assumptions any $(\C,\D)$-triangle is automatically minimal.

\vspace{3mm}
It is not hard to produce subcategories of $\T$ satisfying the orthogonality conditions (TP1').  For instance, start with any $\X \subseteq \T$, and set $\C = {}^{\perp}(\X^{\perp})$ and $\D = \X^{\perp}$.   It is then routine to check that $\D = \C^{\perp}$.  It is an interesting problem to determine which subsets $\X$ of $\T$ will actually produce a torsion pair in this way.  By (the dual of) Proposition 2.3 in \cite{IY}, this is equivalent to asking for which $\X$ is $\X^{\perp}$ covariantly finite in $\T$.  Aihara and Iyama show that this is the case if $\T$ has arbitrary coproducts and $\X$ consists of compact objects \cite{AihIya}, but we cannot assume $\T$ has arbitrary coproducts if we want our theory to apply to the stable category of finitely generated modules over a self-injective algebra.  Similarly, one could take $\C = {}^{\perp}\X$ and $\D = (^{\perp}\X)^{\perp}$ to get another pair of subcategories satisfying (TP1').  We also point out that  both $(^{\perp}\X)^{\perp}$ and $ {}^{\perp}(\X^{\perp})$ are extension-closed subcategories containing $\X$, and hence they both contain $\F(\X)$.  It is natural to wonder when equality might hold.  The following theorem answers both of these questions affirmatively in a special case.

\begin{therm}  Suppose $\X \subseteq \s$ for a simple-minded system $\s$ in $\T$.  Then $({}^{\perp}\X, \F(\X))$ and $(\F(\X), \X^{\perp})$ are torsion pairs in $\T$.  In particular, $\F(\X)$ is a functorially finite subcategory of $\T$.
\end{therm}

\noindent
{\it Proof.}  As we have already observed, the pair $({}^{\perp}\X, (^{\perp}\X)^{\perp})$ satisfies (TP1').  We will check (TP2) for $X \in (\s)_n$ by induction on $n$.  Simultaneously, we will see that $\bfb(X) \in \F(\X)$.  If $n=1$, then we use the triangle $0 \longrightarrow X \stackrel{1}{\longrightarrow} X \rightarrow$ if $X \in \X$, or the triangle $X \stackrel{1}{\longrightarrow} X \longrightarrow 0 \rightarrow$ if $X \in \s \setminus \X$.  Now assume that a $({}^{\perp}\X, (^{\perp}\X)^{\perp})$-triangle exists for every $Y \in (\s)_{n-1}$, and moreover that $\bfb(Y) \in \F(\X)$ for all such $Y$.  If $X \in {}^{\perp}\X$, there is nothing to prove.  Thus we may assume, using Lemma 2.6, that there is a triangle $Y \longrightarrow X \stackrel{f}\longrightarrow S_i \rightarrow $ with $S_i \in \X$ and $Y \in (\s)_{n-1}$.  By the induction hypothesis and the octahedral axiom, we have a commutative diagram where the rows and columns are triangles.
$$\xymatrix{\bfa Y  \ar[d] \ar@{=}[r] & \bfa Y \ar[d] \\ Y \ar[d] \ar[r] & X \ar[d] \ar[r] & S_i \ar[r] \ar@{=}[d] & \\ \bfb Y \ar[d] \ar[r] & W \ar[d] \ar[r] & S_i \ar[r]  & \\ & }$$
The bottom row shows that $W \in \F(\X)$, and the second column yields the desired triangle for $X$.   Moreover, $\bfb X$ is isomorphic to a direct summand of $W$ and thus belongs to $\F(\X)$ by Lemma 2.7.  Finally, if $X \in (^{\perp}\X)^{\perp}$, then $X \cong \bfb X \in \F(\X)$, and it follows that $({}^{\perp}\X, \F(\X)) = ({}^{\perp}\X, (^{\perp}\X)^{\perp})$ is a torsion pair.  The proof that $(\F(\X), \X^{\perp})$ is a torsion pair is entirely analogous. $\Box$ \\

\noindent
{\bf Remark.}  Observe that it follows from the above proof that $\bfa X, \bfb X \in (\s)_n$ whenever $X \in (\s)_n$, for either of the torsion pairs in the theorem.\\

In much of what follows, we will assume that $\T$ has {\it Serre duality}, meaning that there is an exact auto-equivalence $\nu : \T \rightarrow \T$, called a {\it Serre functor}, together with natural isomorphisms $\T(X,Y) \cong D\T(Y, \nu X)$ for all $X, Y \in \T$, where $D = \Hom_k(-,k)$ is the standard duality with respect to the ground field $k$.

\begin{lemma}  Assume $\T$ has a Serre functor $\nu$.  Then, for any $\X \subseteq \T$, we have $\X^{\perp} = \nu ({}^{\perp}\X)$ and $({}^{\perp}\X)^{\perp} = {}^{\perp}(\X^{\perp})$.
\end{lemma}

\noindent
{\it Proof.}  Let $X \in \X$.  Since $\T(Y,X) = 0$ if and only if $\T(X,\nu Y)=0$, we see that $Y \in {}^{\perp}\X$ if and only if $\nu Y \in \X^{\perp}$.  Applying this identity to ${}^{\perp}\X$, we get $({}^{\perp}\X)^{\perp} = \nu({}^{\perp}({}^{\perp}\X)) = {}^{\perp}\nu({}^{\perp}\X) = {}^{\perp}(\X^{\perp})$.  $\Box$ \\

\section{Mutation of simple-minded systems}
\setcounter{equation}{0}

We now focus on simple-minded systems and the associated torsion pairs  in triangulated categories with Serre duality.  Thus, from now on we assume that $\nu$ is a Serre functor for $\T$, that $\s \subseteq \T$ is a simple-minded system with $\nu(\s[1]) = \s$ and that $\X \subseteq \s$ also satisfies $\nu(\X[1]) = \X$.  From Theorem 3.3, we have two torsion pairs $({}^{\perp}\X, \F(\X))$ and $(\F(\X),\X^{\perp})$, and we define $\bfa : \T \rightarrow {}^{\perp}\X$, $\bfb, \bfc : \T \rightarrow \F(\X)$ and $\bfd : \T \rightarrow \X^{\perp}$ via the minimal triangles $$\bfa X \rightarrow X \rightarrow \bfb X \rightarrow\ \ \ \mbox{and}\ \ \ \bfc X \rightarrow X \rightarrow \bfd X \rightarrow$$ corresponding to these two torsion pairs respectively.  Additionally, we set $\alpha = \bfa \circ [-1]$ and $\beta = \bfd \circ [1]$.  We begin by defining mutations of simple-minded systems, generalizing our previous definition from \cite{MutSym}.

\begin{defin} With the notation and assumptions above, we define the {\bf left mutation of $\s$ at $\X$} to be $$\mu^+_{\X}(\s) = \X \cup \{\alpha S_i\ |\ S_i \in \s \setminus \X\}.$$  
Dually, we define the {\bf right mutation of $\s$ at $\X$} to be $$\mu^-_{\X}(\s) = \X \cup \{\beta S_i\ |\ S_i \in \s \setminus \X\}.$$
For $S_i \in \s$, we may also use $\mu^+_{\X}(S_i)$ to denote either $S_i$ or $\alpha S_i$, respectively, depending on whether or not $S_i \in \X$, and similarly for $\mu^-_{\X}(S_i)$.
\end{defin}

\noindent
{\bf Remark.}  In \cite{KoYa} (see also \cite {KeYa}), Koenig and Yang define mutation for {\it simple-minded collections}, which can be viewed as the derived category analogues of simple-minded systems.  Their definition suggests a mutation of simple-minded systems that differs from ours by an application of the suspension functor:  namely, $\mu^+_{\X}(S_i) = S_i[1]$ if $S_i \in \X$, while $\mu^+_{\X}(S_j)$ is defined to be the cone of the left $\F(\X)$ approximation of $S_j[-1]$ if $S_j \in \s \setminus \X$ (see (5.3)).  In analogy with the mutation of quivers, this latter definition is more logical: mutating a quiver at a single vertex is a local operation, affecting the arrows in a neighborhood of that vertex and leaving the farther reaches of the quiver unchanged.  Likewise this definition of mutation of simple-minded systems replaces the ``simples'' in $\X$ with their suspensions and will have no affect on the ``simples'' that are ``sufficiently far'' from those in $\X$.   However, we have opted for our alternate definition here as it is more convenient for computing examples: especially, if $\X$ is large, in which case we only need to compute the affect of mutating on the objects outside of $\X$.  Moreover, we can apply our mutation procedure repeatedly with respect to the same subset $\X$ using a single torsion pair.\\

The following theorem represents the main result of this section.  It answers a question raised in \cite{MutSym} and provides a useful way of obtaining many nontrivial examples of simple-minded systems in stable categories.

\begin{therm} Suppose $\s$ is a simple-minded system in $\T$ and $\X \subseteq \s$ satisfies $\nu(\X[1]) = \X$.  Then $\mu^+_{\X}(\s)$ and $\mu^-_{\X}(\s)$ are simple-minded systems in $\T$.
\end{therm}

The proof will be completed in several steps.  We begin with some simple observations.

\begin{lemma}  Assume $\X \subseteq \s$ satisfies $\nu(\X[1]) = \X$.  Then $\nu(\F(\X)[1]) = \F(\X)$.  Furthermore, $\nu(\bfa X)[1] \cong \bfa (\nu X[1])$ and $\nu(\bfb X)[1] \cong \bfb (\nu X[1])$ for all $X \in \T$, and similarly for $\bfc$ and $\bfd$.
\end{lemma}

\noindent
{\it Proof.}  We show that $\F(\X) \subseteq \nu(\F(\X)[1])$, as the proof of the reverse inclusion is similar.  For $X \in (\X)_n$, we will prove that $X \in \nu(\F(\X)[1])$ by induction on $n$.  If $n = 1$, then $X \in \X = \nu(\X[1]) \subseteq \nu(\F(\X)[1])$.  Now assume that $(\X)_{n-1} \subseteq \nu(\F(\X)[1])$ for some $n \geq 2$.  If $X \in (\X)_n$, we have a triangle $X' \rightarrow X \rightarrow S_i \rightarrow$ with $X' \in (\X)_{n-1} \subseteq \nu(\F(\X)[1])$ and $S_i \in \X = \nu(\X[1])$.  Thus, we also have a triangle $\nu^{-1}(X')[-1] \rightarrow \nu^{-1}(X)[-1] \rightarrow \nu^{-1}(S_i)[-1] \rightarrow$, which shows that $\nu^{-1}(X)[-1] \in \F(\X)$, and thus that $X \in \nu(\F(\X)[1])$.   
 
 Now consider a triangle $\bfa X \rightarrow X \rightarrow \bfb X \rightarrow$, to which we can apply the exact functor $\nu \circ -[1]$.  Noting that $\nu((\bfa X)[1]) \in {}^{\perp}\nu(\F(\X)[1]) = {}^{\perp}\F(\X)$, we see that the resulting triangle is isomorphic to $\bfa(\nu X[1]) \rightarrow \nu X[1] \rightarrow \bfb(\nu X[1]) \rightarrow$. $\Box$ \\

\begin{lemma}  For any $M \in \ppX$, we have $\beta \alpha M \cong M \cong \alpha \beta M$.
\end{lemma}

\noindent
{\it Proof.}  We may assume that $M$ is indecomposable.  By definition of $\alpha$, we have a minimal triangle $\alpha M \stackrel{f}{\rightarrow} M[-1] \rightarrow \bfb(M[-1]) \rightarrow$.  Notice that $f \neq 0$ since $M [-1] \in {}^{\perp}(\X[-1]) = {}^{\perp}(\nu \X) = \X^{\perp}$ implies that $M[-1] \notin \F(\X)$.  Rotating this triangle to the right twice yields another triangle $\bfb ( M[-1]) \rightarrow (\alpha M)[1] \stackrel{f[1]}{\rightarrow} M \rightarrow$, which is again minimal since $M$ is indecomposable and $f[1] \neq 0$.    Since $\bfb (M[-1]) \in \F(\X)$ and $M \in \X^{\perp}$, we must have $M \cong \bfd ((\alpha M)[1]) = \beta \alpha M$.  The other identity is proved similarly.  $\Box$ \\

Since any $S_i \in \s \setminus \X$ belongs to $\ppX$, we obtain the following corollary, which shows that left and right mutation are mutually inverse operations.  We point out that the operation $\mu^-_{\X}$ (resp. $\mu^+_{\X}$) is still defined on $\mu^+_{\X}(\s)$ (resp. on $\mu^-_{\X}(\s)$) even without knowing that the latter is a simple-minded system.

\begin{coro} Let $\s$ be a simple-minded system and suppose that $\X \subseteq \s$ satisfies $\nu(\X[1]) = \X$.  Then $$\mu^-_{\X}(\mu^+_{\X}(\s)) = \s = \mu^+_{\X}(\mu^-_{\X}(\s)).$$
\end{coro}

\begin{lemma} For any triangle $\bfa M \stackrel{f}{\rightarrow} M \stackrel{g}{\rightarrow} \bfb M \rightarrow$ and any $X \in \X$, we have
\begin{enumerate}
\item[(a)] The map $\T(g,X) : \T(\bfb M, X) \rightarrow \T(M,X)$ is an isomorphism;
\item[(b)] The map $\T(X,f) : \T(X, \bfa M) \rightarrow \T(X,M)$ is a monomorphism;
\item[(c)] If $M \in \X^{\perp}$, then $\bfa M \in \ppX$.
\end{enumerate}
\end{lemma}

\noindent
{\it Proof.}  (a) Since $g$ is a left $\F(\X)$-approximation, we know that $\T(g,X)$ is onto.   To see that it is one-to-one, suppose we have a map $h : \bfb M \rightarrow X$ such that $hg = 0$.   By the octahedral axiom we have a commutative diagram where the rows and columns are triangles.
$$\xymatrix{& M \ar[d]^g \ar@{=}[r] & M \ar[d]^{0} \\ Y \ar@{=}[d] \ar[r] & \bfb M \ar[d] \ar[r]^{h} & X \ar[r] \ar[d] & \\ Y  \ar[r] & (\bfa M)[1] \ar[d] \ar[r] & X \oplus M[1] \ar[r]  \ar[d] & \\ & & }$$
Twice rotating the bottom triangle to the left yields a triangle $\bfa M \rightarrow X[-1] \oplus M \rightarrow Y \rightarrow$ where $Y \in \F(\X)$ by Lemma 2.6.  Since this triangle must be isomorphic to the direct sum of the minimal $({}^{\perp}\X, \F(\X))$-triangles for $X[-1]$ and $M$, we see that the component of the second map $X[-1] \rightarrow Y$ must be a split monomorphism and $X[-1] \in \F(\X)$.   Thus the triangle $X[-1] \rightarrow Y \rightarrow \bfb M \rightarrow X$, obtained by rotating the triangle in the first row above, splits and it follows that $h=0$.

(b) Since $\T(g,X)$ is injective, $\T(f[1],X)$ must be surjective.  We fix an isomorphism $h : \bfa (\nu M[1]) \rightarrow \nu (\bfa M)[1]$ by Lemma 4.3, and we obtain a commutative diagram 

$$\xymatrixcolsep{4.0pc} \xymatrix{ \T(M[1],X) \ar[r]^{\T(f[1],X)} \ar[d]^{\nu}_{\cong} & \T( (\bfa M)[1], X) \ar[d]^{\nu}_{\cong} \\ 
\T(\nu M[1], \nu X) \ar[r]^{\T(\nu f[1],\nu X)} \ar@{=}[d] & \T( \nu (\bfa M)[1], \nu X) \ar[d]^{\T(h,\nu X)}_{\cong} \\ 
\T(\nu M[1], \nu X) \ar[r]^{\T(\nu f[1] \circ h,\nu X)} \ar[d]^{}_{\cong} & \T( \bfa (\nu M[1]), \nu X) \ar[d]_{\cong} \\ D\T(X, \nu M[1]) \ar[r]^{D\T(X,\nu f[1] \circ h)} & D\T(X, \bfa (\nu M[1])).
}$$
Since the morphism on the bottom is surjective, the map $\T(X, \nu f[1] \circ h)$ must be injective.  This establishes (b) for the object $\nu M[1]$ (in place of $M$), but every object of $\T$ arises in this way.

(c) This is an immediate consequence of (b).  $\Box$ \\

For convenience, we state the dual of Lemma 4.6 without proof.

\begin{lemma}  For any triangle $\bfc M \stackrel{f}{\rightarrow} M \stackrel{g}{\rightarrow} \bfd M \rightarrow$ and any $X \in \X$, we have
\begin{enumerate}
\item[(a)] The map $\T(X,f) : \T(X, \bfc M) \rightarrow \T(X,M)$ is an isomorphism;
\item[(b)] The map $\T(g,X) : \T( \bfd M, X) \rightarrow \T(M,X)$ is a monomorphism;
\item[(c)] If $M \in {}^{\perp}\X$, then $\bfd M \in \ppX$.
\end{enumerate}
\end{lemma}

Observe that, as a consequence of the two previous lemmas along with Lemma 3.4, $\alpha :  {}^{\perp}\X \rightarrow \ppX$ and $\beta : \X^{\perp} \rightarrow \ppX$.  In fact, $\alpha$ and $\beta$ are functors when restricted to these subcategories.  To see this, it suffices to check that $\bfa : \X^{\perp} \rightarrow \ppX$ and $\bfd : {}^{\perp}\X \rightarrow \ppX$ are well-defined on morphisms.  For $\bfa$, this follows from the fact that
$\T(\bfa M, (\bfb N)[-1]) \cong D\T(\nu^{-1} (\bfb N)[-1], \bfa M) \cong D\T(\bfb(\nu^{-1}N[-1]), \bfa M) = 0$, as $\bfa M \in \ppX$ for any $M \in \X^{\perp}$.  The argument for $\bfd$ is similar.

\begin{propos}  Suppose $\s$ is a simple-minded system with $\X \subseteq \s$ satisfying $\nu(\X[1]) = \X$.  Then $\mu^+_{\X}(\s)$ and $\mu^-_{\X}(\s)$ are sets of orthogonal bricks in $\T$.
\end{propos}

\noindent
{\it Proof.}  We only provide the proof for $\mu^+_{\X}(\s)$, as the proof for $\mu^-_{\X}(\s)$ is similar.  First notice that $\alpha S_i \in \ppX$ for all $S_i \in \s \setminus \X$ by the above remarks.  The rest of our proof is modeled on the proof of Theorem 6.2 in \cite{MutSym}.   

Suppose that $S_i, S_j \in \s \setminus \X$, and consider the exact diagram obtained by applying $\T(\alpha S_j,-)$ and $\T(-,S_i[-1])$ to the triangles defining $\alpha S_i$ and $\alpha S_j$ respectively.
$$\xymatrix{ & & \T(\bfb (S_j[-1])[-1], S_i[-1]) \\ \T(\alpha S_j, \bfb (S_i[-1])[-1]) \ar[r]^-0 &\T(\alpha S_j, \alpha S_i) \ar[r] & \T(\alpha S_j, S_i[-1]) \ar[r]^0 \ar[u]^0 & \T(\alpha S_j, \bfb(S_i[-1])) \\ & & \T(S_j[-1], S_i[-1]) \ar[u] \\ & & \T(\bfb(S_j[-1]), S_i[-1]) \ar[u]^0 }$$
To get the zero maps we have used that $\T(\alpha S_j, \bfb (S_i[-1])[-1]) \cong \T(\alpha S_j, \nu \bfb (\nu^{-1}S_i[-2])) \cong D\T(\bfb(\nu^{-1}S_i[-2]), \alpha S_j) = 0$ and $\T(\alpha S_j,\bfb(S_i[-1]))=0$ since $\alpha S_j \in \ppX$; that $\T(\bfb(S_j[-1])[-1], S_i[-1]) \cong \T(\bfb(S_j[-1]), S_i) =0 $ since $S_i \in \ppX$;  and that $\T(\bfb(S_j[-1]),  S_i[-1]) \cong D\T(\nu^{-1}(S_i[-1]), \bfb(S_j[-1])) = 0$ since $\nu^{-1}(S_i[-1]) \in \nu^{-1}(\ppX[-1]) = \ppX$.  Thus the remaining nonzero maps now yield isomorphisms 
$$\T(\alpha S_j, \alpha S_i) \cong \T(\alpha S_j, S_i[-1]) \cong \T(S_j[-1], S_i[-1]) \cong \T(S_j,S_i),$$ which is a division ring if $i=j$ and vanishes otherwise.  $\Box$ \\

\vspace{2mm}
We need two additional lemmas in order to prove Theorem 4.2.

\begin{lemma} For any $N \in {}^{\perp}\X$, $\alpha N$ is a direct summand of $\alpha \bfd N$.  
\end{lemma}

\noindent
{\it Proof.}  By definition of $\alpha$ we have a triangle $\alpha N \rightarrow N[-1] \rightarrow \bfb (N[-1]) \rightarrow $.  Rotating to the right now yields the triangle in the top row below.
$$\xymatrix{\bfb(N[-1]) \ar[r] \ar@{-->}[d] & (\alpha N) [1] \ar@{=}[d] \ar[r] & N \ar@{-->}[d]^{\varphi} \ar[r] & \bfb(N[-1])[1] \ar[r] \ar@{-->}[d] & \\ \bfc((\alpha N)[1]) \ar[r]_g & (\alpha N)[1] \ar[r] & \beta \alpha N \ar[r] & \bfc((\alpha N)[1])[1] \ar[r] & } $$
The triangle in the second row is that used to define $\beta \alpha N$, and the morphism of triangles exists since $\bfb(N[-1]) \in \F(\X)$ while $g$ is a right $\F(\X)$-approximation.  Applying the inverse of the suspension functor to the diagram now yields the morphism of triangles. 
$$\xymatrix{\alpha N \ar@{=}[d] \ar[r] & N[-1] \ar[d]^{\varphi[-1]} \ar[r] & \bfb(N[-1]) \ar[r] \ar[d] & \\  \alpha N \ar[r] & (\beta \alpha N)[-1] \ar[r] & \bfc( \alpha N[1]) \ar[r] & } $$
Since $\alpha N \in {}^{\perp}\X$ while $\bfb(N[-1]), \bfc(\alpha N[1]) \in \F(\X)$, these are both $({}^{\perp}\X, \F(\X))$-triangles.  Moreover, they are minimal since $\alpha N \in \ppX$, which forces $(\alpha N)[1] \in \ppX [1] = \nu^{-1}(\ppX) \subseteq {}^{\perp}\X$. Consequently, we see that $\alpha(\varphi) = \bfa (\varphi[-1])$ is an isomorphism.  Now let $f : N \rightarrow \bfd N$ be a minimal left $\X^{\perp}$-approximation.  Since $\beta \alpha N \in \X^{\perp}$, $\varphi$ must factor through $f$.  As $\bfd N$ also belongs to $\ppX \subseteq {}^{\perp}\X$, the isomorphism $\alpha(\varphi)$ factors through $\alpha(f)$, which forces $\alpha(f)$ to be a split monomorphism.  $\Box$ \\

\begin{lemma}  Suppose that $L \stackrel{f}{\rightarrow} M \stackrel{g}{\rightarrow} N \stackrel{h}{\rightarrow}$ is a triangle with $N \in \F(\X)$ and $L \in \X^{\perp}$.  Then $\bfa L \cong \bfa M$.
\end{lemma}

\noindent
{\it Proof.}  Let $i_L : \bfa L \rightarrow L$ and $i_M : \bfa M \rightarrow M$ be the natural maps.  Since $N \in \F(\X)$, we have $gi_M = 0$ and thus $i_M$ factors through $f$.  Since $\bfa M \in {}^{\perp}\X$, we see that $i_M$ factors through $ f i_L$.  Thus $f i_L$ is a right ${}^{\perp}\X$-approximation of $M$.  We claim that $f i_L$ is right minimal.  Otherwise, there is a nonzero direct summand $L' $ of $\bfa L$ for which $f i_L|_{L'} = 0$, which implies that $i_L|_{L'}$ factors through the map $h[-1] : N[-1] \rightarrow L$.  However, $\T(L', N[-1]) \cong D\T(\nu^{-1}N[-1], L')= 0$ as $\bfa L$, and hence $L'$, belongs to $\ppX$ by Lemma 4.6(c) and  $\nu^{-1}N[-1] \in \F(\X)$ by Lemma 4.3.  This implies that $i_L|_{L'} = 0$, contradicting the minimality of $i_L$.  Thus $f i_L$ is a minimal right ${}^{\perp}\X$-approximation of $M$, and we have $\bfa M \cong \bfa L$.  $\Box$ \\

\noindent
{\it Proof of Theorem 4.2.}  In light of Proposition 4.8, it only remains to show that $\F(\mu^+_{\X}(\s)) = \T$.  We begin by showing that $\alpha(\ppX) \subseteq \F(\mu^+_{\X}(\s))$ by induction on the $\s$-length of an object in $\ppX$.    If $M \in \ppX \cap (\s)_1$, then $M \cong S_i \in \s \setminus \X$ and $\alpha M \cong \alpha S_i \in (\mu^+_{\X}(\s))_1$.  Now suppose that $\alpha N \in \F(\mu^+_{\X}(\s))$ for all $N \in \ppX \cap (\s)_{n-1}$, and let $M \in \ppX \cap (\s)_n$.  We can find a triangle $S_i \stackrel{f}{\rightarrow} M \stackrel{g}{\rightarrow} N \rightarrow$ with $S_i \in \s \setminus \X$, $f \neq 0$ and $N \in (\s)_{n-1}$.  Since $M, S_i[1] \in {}^{\perp} \X$, we must have $N \in {}^{\perp}\X$ as well.  In particular, Lemma 4.9 shows that $\alpha N \in \F(\mu^+_X(\s))$ since $\bfd N \in \ppX \cap (\s)_{n-1}$ and $\F(\mu^+_X(\s))$ is closed under direct summands thanks to Proposition 4.8 and Lemma 2.7.

If $\alpha$ were an exact functor on ${}^{\perp}\X$, we would now be done.  Unfortunately, $\alpha$ is typically not exact.  Instead, we construct a triangle  $L \rightarrow \alpha M \rightarrow \alpha N \rightarrow$ with $L \in \F(\mu^+_{\X}(\s))$ as follows.  First, we use the octahedral axiom to obtain a commutative diagram in which the rows and columns are triangles.
$$\xymatrix{ S_i[-1] \ar@{=}[d] \ar[r] & M' \ar[d]^{} \ar[r]^{}  & \alpha N \ar[r] \ar[d] & \\ S_i[-1]  \ar[r]^{f[-1]} &  M[-1] \ar[d] \ar[r]^{g[-1]} & N[-1] \ar[r]  \ar[d] & \\ & \bfb(N[-1]) \ar[d] \ar@{=}[r] & \bfb(N[-1]) \ar[d] \\ & & }$$
Since $N, S_i \in {}^{\perp}\X$, we know that $\alpha N \in \ppX$ and $S_i[-1] \in \X^{\perp}$, and thus the triangle in the top row shows that $M' \in \X^{\perp}$.   By Lemma 4.10, we have $\bfa M' \cong \bfa(M[-1]) = \alpha M$.  Applying the octahedral axiom again we have the commutative diagram
$$\xymatrix{ L \ar[d] \ar[r] & \alpha M \ar[d] \ar[r] & \alpha N \ar[r] \ar@{=}[d] & \\ S_i[-1] \ar[d] \ar[r] & M' \ar[d] \ar[r] & \alpha N \ar[r]  & \\ \bfb M' \ar[d] \ar@{=}[r] & \bfb M' \ar[d] \\ & }$$
which yields the desired triangle $L \rightarrow \alpha M \rightarrow \alpha N \rightarrow$ as its top row.  Moreover, this triangle shows that $L \in \X^{\perp}$ since $\alpha M, \alpha N \in \ppX$ and thus $(\alpha N) [-1] \in \ppX[-1] = \nu(\ppX) \subseteq \X^{\perp}$.  Now applying Lemma 4.10 to the triangle in the first column shows that $\bfa L \cong \bfa(S_i[-1]) = \alpha S_i$.  In particular, we have a triangle $\alpha S_i \rightarrow L \rightarrow \bfb L \rightarrow$, showing that $L \in \F(\mu^{+}_{\X}(\s))$.  It follows that $\alpha M \in \F(\mu^+_{\X}(\s))$, completing the induction argument.

Finally, we know that $\T = {}^{\perp}\X * \F(\X)$ with ${}^{\perp}\X \subseteq \F(\X) * \ppX$ by Lemma 4.7, and it follows from Lemma 4.4 and the remarks following Lemma 4.7 that $\ppX = \alpha \beta (\ppX) \subseteq \alpha (\ppX)$.  Thus $\ppX$ and $\F(\X)$ are contained in $\F(\mu^+_{\X}(\s))$, and hence $\T$ is as well.  $\Box$ \\

\section{Okuyama complexes}

We now set out to show that mutation of simple-minded systems inside stable module categories corresponds to the derived equivalences induced by Okuyama complexes, giving an alternative proof of the unpublished lemma of Okuyama cited in the introduction.   We start by fixing our notation.  Let $\Lambda$ be a basic self-injective algebra over an algebraically closed field $k$, and set $\s = \{S_1, \ldots, S_n\}$, a complete set of representatives of the isomorphism classes of the simple (right) $\Lambda$-modules.   We also write $P_i = e_i \Lambda$ for the indecomposable projective and primitive idempotent corresponding to $S_i$, and for any subset $U \subseteq \s$ we will write $e_U = \sum_{S_i \in U} e_i$, $P_U = e_U \Lambda$ and $Q_U = (1-e_U)\Lambda$.  We let $\N = D\Hom_{\Lambda}(-,\Lambda) \cong -\otimes_{\Lambda} D\Lambda$ denote the Nakayama functor, which induces an auto equivalence of $\stmod \Lambda$ and satisfies $\N \cong \nu \Omega^{-1}$, where $\nu$ is the Serre functor for $\stmod \Lambda$ and $\Omega$ is the syzygy functor.  Thus we can mutate $\s$ with respect to any subset $\X \subset \s$ that is invariant under $\N$.  Furthermore, corresponding to any such subset of $\s$ we have an Okuyama tilting complex $T_{\X} \in K^b(\proj \Lambda)$ given by $T_{\X} = \oplus_{i=1}^n T_i$ where
\begin{eqnarray} T_i & = & \left\{ \begin{array}{rl} P_i[-1], &  \mbox{if\ } i \notin \X \\ \ul{P_i} \rightarrow L_i, & \mbox{if\ } i \in \X
\end{array} \right.
\end{eqnarray}
with the maps $P_i \rightarrow L_i$ being minimal left $\add(Q_\X)$-approximations, and with $P_i$ underlined to signify that it is the degree $0$ term of this complex.  The assumption that $\N(\X) = \X$ ensures that $T_{\X}$ is a tilting complex.  We thus set $\Gamma = \Endo_{K^b(\Lambda)}(T_{\X})$, which is called the {\it tilting mutation} $\mu^+_{\X}(\Lambda)$ in \cite{MutSym}.  We write $F : D^b(\rmod \Gamma) \rightarrow D^b(\rmod \Lambda)$ for the induced derived equivalence and $\ul{F} : \stmod \Gamma \rightarrow \stmod \Lambda$ for the induced stable equivalence as given by \cite{DCSE, HuXi3}.  %We note that $\ul{F}$ is really only well-defined up to a power of the suspension functor in general
We can now state the main result of this section.\\

\noindent
{\bf Okuyama's Lemma.}  [cf. Lemma 2.1' in \cite{Oku}]  \emph{With the notation and assumptions above, $\ul{F}$ takes the set of simple $\Gamma$-modules to $\mu^+_{\X}(\s)$ in $\stmod \Lambda$.}\\

The proof we outline here relies mainly on recent work of Koenig and Yang \cite{KoYa} connecting mutation of simple-minded collections in $D^b(\rmod \Lambda)$ to the mutation of silting objects (recently studied by Aihara and Iyama \cite{AihIya}) in $K^b(\proj \Lambda)$.  Thus we begin by reviewing the relevant definitions and results that we will need.  We will continue to write $\N$ for the auto-equivalences induced by the Nakayama functor on $K^b(\proj \Lambda)$ and on $D^b(\rmod \Lambda)$.  Thus in $D^b(\rmod \Lambda)$, $\N = -\otimes_{\Lambda}^{\mathbb{L}} D\Lambda$, and our usage of $\N$ will be clear from context.  We note that $\N$ is a Serre functor for $K^b(\proj \Lambda)$.

We say that an object $M$ in a triangulated category $\T$ is a {\it silting object} if $\add(M)$ generates $\T$ (as a triangulated category) and $\T(M,M[i]) = 0$ for all $i>0$.  Aihara and Iyama introduced a mutation procedure for silting objects in \cite{AihIya}.  Namely, if $M = U \oplus V$ is a silting object in $\T$, then one obtains a new silting object $\mu^+_U(M)$ by replacing $V$ with the cone $N_V$ of a left $\add(U)$-approximation $V \stackrel{f}{\rightarrow} D$ in $\T$.  That is, $\mu^+_U(M) = U \oplus N_V$.  In case $\T = K^b(\proj \Lambda)$ for a self-injective algebra $\Lambda$, it is not difficult to see that any $\N$-stable silting object $M$ (i.e., satisfying $\N(M) \cong M$) must be a tilting complex.  Furthermore, if $\N(U) = U$ for a summand $U$ of $M$, then the silting mutation $\mu^+_U(M)$ remains $\N$-stable and thus a tilting complex.  In fact, as shown in Section 2.7 of \cite{AihIya}, Okuyama tilting complexes are precisely those obtained by performing $\N$-stable silting mutations to the stalk complex $\Lambda \in K^b(\proj \Lambda)$.   In our setting, one easily checks that the complex $T_{\X}$ described above is realized by mutating the tilting complex $\Lambda$ at $P_{\X}$ and de-suspending: 
  \begin{eqnarray} T_{\X} \cong \mu^+_{P_{\X}}(\Lambda)[-1]. \end{eqnarray}

In \cite{KoYa}, Koenig and Yang define a mutation procedure for simple-minded collections (which were essentially introduced by Rickard in \cite{EDCSA}), and show that it is compatible with mutation of silting objects in an appropriate sense.  We say that a set of objects $\{X_1, \ldots, X_r\} \subseteq \T$ is a {\it simple-minded collection} if 
\begin{itemize}
\item $\T(X_i, X_j[m])  = 0$ for all $m <0$; 
\item $\T(X_i, X_j) = 0$ for all $i \neq j$ and $\T(X_i,X_i) \cong k$ for all $i$; and 
\item $X_1, \ldots, X_n$ generate $\T$ as a triangulated category.  
\end{itemize}
Clearly the set of simple $\Gamma$-modules $\s' = \{S'_1, \ldots, S'_n\}$ forms a simple-minded collection in $D^b(\rmod \Gamma)$.  Moreover, this simple-minded collection is $\N$-stable.  Since $F : D^b(\rmod \Gamma) \rightarrow D^b(\rmod \Lambda)$ is an equivalence, we see that $F(\s') = \{F(S'_1), \ldots, F(S'_n)\}$ is a  simple-minded collection in $D^b(\rmod \Lambda)$.  In fact, the latter simple-minded collection must also be $\N$-stable.  To see this, we can replace $F$ with a standard derived equivalence $F'$ (i.e., induced by tensoring with a two-sided tilting complex)  that agrees with $F$ on isomorphism classes of objects by Corollary 3.5 of \cite{DEDF} and then apply Proposition 5.2 in \cite{DEDF} to $F'$ to see that $F$ must commute with $\N$ at least on isomorphism classes of objects. 

Koenig and Yang have generalized this procedure to define a bijection $\phi_{21}$ from the set of equivalence classes of silting objects in $K^b(\proj \Lambda)$ to the set of equivalence classes of simple-minded collections in $D^b(\rmod \Lambda)$  (here, two (sets of) objects are said to be `equivalent' if they generate the same additive subcategory).   Furthermore, in our setting where $\Lambda$ is self-injective, Al-Nofayee has previously shown that this construction of a silting object corresponding to an $\N$-stable simple-minded collection produces a tilting complex in $K^b(\proj \Lambda)$ \cite{AlNo},  i.e., an  $\N$-stable silting complex.  Thus their map $\phi_{21}$ restricts to a bijection between the set of equivalence classes of tilting complexes in $K^b(\proj \Lambda)$ and the set of equivalence classes of $\N$-stable simple-minded collections in $D^b(\rmod \Lambda)$.  By construction this bijection commutes with the suspension functors of these two categories; i.e., $\phi_{21}(T[1]) \cong \phi_{21}(T)[1]$ for any tilting complex $T \in K^b(\proj \Lambda)$.  The definition and main properties of the map $\phi_{21}$ that we need will be summarized in the theorem below.

While mutation for simple-minded collections may not work as generally as mutation for simple-minded systems, Koenig and Yang prove that the former mutations are always possible in the bounded derived category of a finite-dimensional algebra.  While they prove this for mutations with respect to a single object of the simple-minded collection, the same definition and argument carry over to the setting of mutations with respect to any subset of the collection.   Thus, modifying \cite{KoYa}, for a simple-minded collection $X_1, \ldots, X_n \in D^b(\rmod \Lambda)$ and $\X \subseteq \{X_1, \ldots, X_n\}$ we obtain another simple-minded collection consisting of the objects $\mu^+_{\X}(X_i)$ defined as follows:  
\begin{eqnarray} \mu^+_{\X}(X_i) = \left\{ \begin{array}{rl} X_i[1], & \mbox{if}\ X_i \in \X \\ cone(X_i[-1] \stackrel{g_i}{\rightarrow} Y_i), & \mbox{if}\ X_i \notin \X
\end{array} \right. \end{eqnarray} 
where the map $g_i$ is a minimal left $\F(\X)$-approximation of $X_i[-1]$, and $\F(\X)$ is the extension-closure of $\X$ in $\T = D^b(\rmod \Lambda)$.   (Although we have used the same notation for mutations of simple-minded collections as we have previously for simple-minded systems, the context will determine which one is meant, as we only consider simple-minded collections in $D^b(\rmod \Lambda)$ and simple-minded systems in $\stmod \Lambda$.)

We can now summarize the results that we will use from \cite{KoYa}, with part (4) corresponding to part of Theorem 7.12. %8.11.

\begin{therm}[Koenig, Yang \cite{KoYa}]  Let $T \in K^b(\proj \Lambda)$ be a basic tilting complex and $F : D^b(\rmod \Endo(T)) \rightarrow D^b(\rmod \Lambda)$ the corresponding derived equivalence.  Then the image of the set of simple $\Endo(T)$-modules under $F$ is an $\N$-stable simple-minded collection in $D^b(\rmod \Lambda)$, denoted $\phi_{21}(T)$.  The map $\phi_{21}$ induces a bijection between the set of equivalence classes of basic tilting complexes in $K^b(\proj \Lambda)$ and the set of equivalence classes of $\N$-stable simple-minded collections in $D^b(\rmod \Lambda)$, and satisfies:
\begin{enumerate}
\item $\phi_{21}(\Lambda) = \s$, the set of simple $\Lambda$-modules;
\item $\phi_{21}(T[1]) = \phi_{21}(T)[1]$; 
\item There is a natural bijection between the indecomposable summands of $T$ and the objects in $\phi_{21}(T)$; 
\item For any direct summand $U$ of $T$ and the corresponding subset $\X$ of $\phi_{21}(T)$, we have $$\phi_{21}(\mu^+_{U}(T)) = \mu^+_{\X}(\phi_{21}(T)).$$ 
\end{enumerate}
\end{therm}

\begin{coro} Let $\s = \{S_1, \ldots, S_n\} \subset \rmod \Lambda \subset D^b(\rmod \Lambda)$ be the simple-minded collection consisting of the set of simple $\Lambda$-modules.  Assume $\X \subseteq \s$ satisfies $\N(\X) = \X$ and let $T_{\X}$ be the corresponding Okuyama tilting complex as in (5.1).  Then the induced derived equivalence $F : D^b(\rmod \Gamma) \rightarrow D^b(\rmod \Lambda)$ sends the simple $\Gamma$-modules to $\mu^+_{\X}(\s)[-1]$.
\end{coro}

\noindent
{\it Proof.}  We have $$F(\s') = \phi_{21}(T_{\X}) = \phi_{21}(\mu^+_{P_{\X}}(\Lambda)[-1]) = \mu^+_{\X}(\phi_{21}(\Lambda))[-1] = \mu^+_{\X}(\s)[-1]. \ \  \Box$$

Finally, we show that mutating the simples in $D^b(\rmod \Lambda)$ lifts the corresponding mutation procedure in $\stmod \Lambda$ as defined in the previous section.  We let $Q_{\Lambda} : D^b(\rmod \Lambda) \rightarrow \stmod \Lambda$ denote the localization functor, which we write as $Q$ when there is no ambiguity.  It is exact and dense, and its restriction to $\rmod \Lambda$ is full.

\begin{propos}  For any $\N$-stable subset $\X$ of the set $\s$ of simple $\Lambda$-modules in $D^b(\rmod \Lambda)$, we have $Q(\mu^+_{\X}(\s)) = \mu^+_{Q(\X)}(Q(\s))[1]$.  \end{propos}

\noindent
{\it Proof.}  First, if $S_i \in \X$, then $Q(\mu^+_{\X}(S_i)) = Q(S_i[1]) = \mu^+_{Q(\X)}(Q(S_i))[1]$ by definition.   So we now assume that $S_i \in \s \setminus \X$.  By definition of $\mu^+_{\X}(\s)$,  $D^b(\rmod \Lambda)$ contains a triangle
 $$\mu^+_{\X}(S_i)[-1] \longrightarrow S_i[-1] \stackrel{g_i}{\longrightarrow} Y_i \longrightarrow \mu^+_{\X}(S_i)$$ where $g_i$ is a minimal left $\F(\X)$-approximation.  Since $Q$ is exact, we also have a triangle 
 \begin{eqnarray} Q(\mu^+_{\X}(S_i))[-1]  \longrightarrow Q(S_i[-1]) \stackrel{Q(g_i)}{\longrightarrow} Q(Y_i) \longrightarrow Q(\mu^+_{\X}(S_i)) \end{eqnarray}
in $\stmod \Lambda$.  It thus remains to show that $Q(g_i)$ is a minimal left $\F(Q(\X))$-approximation.

We first show that $Q(Y_i) \in \F(Q(\X))$.  Since $\X \subset \s \subset \rmod \Lambda$, and the latter is extension-closed in $D^b(\rmod \Lambda)$, the extension closure $\F(\X)$ of $\X$ will be contained in $\rmod \Lambda$ as well.   In fact, it will consist of all the $\Lambda$-modules whose simple composition factors all belong to $\X$.  Since $Q$ takes short exact sequences in $\rmod \Lambda$ to distinguished triangles in $\stmod \Lambda$, it is clear that $Q(\F(\X)) \subseteq \F(Q(\X))$.  Moreover, equality must hold since any distinguished triangle $A \rightarrow B \stackrel{f}{\rightarrow} S \rightarrow$ in $\stmod \Lambda$ with $S$ simple and $f \neq 0$ can be lifted to a short exact sequence $\ses{A}{B}{S}{}{f}$ in $\rmod \Lambda$.

Now suppose that $h : Q(S_i[-1]) \rightarrow M$ with $M \in \F(Q(\X))$.  By the above remarks, we may write $M = Q(M')$ for some $M' \in \F(\X)$.  In order to lift $h$ to a map in the derived category, we will replace $Q(S_i[-1])$ by $Q(\Omega S_i)$ as follows.  The natural map $\varphi : S_i[-1] \rightarrow \Omega S_i$ in $D^b(\rmod \Lambda)$ corresponding to the truncation of the degree-$1$ term of the projective resolution of $S_i[-1]$ induces an isomorphism $Q(\varphi): Q(S_i[-1]) \rightarrow Q(\Omega S_i)$ in $\stmod \Lambda$.  Thus we may find a map $f : \Omega S_i \rightarrow M'$ such that $Q(f) = h Q(\varphi)^{-1}$.  Since $g_i$ is a left $\F(\X)$-approximation, there exists a map $u : Y_i \rightarrow M'$ such that $f \varphi = u g_i$.  Hence, applying $Q$ now yields $h = Q(f) Q(\varphi) = Q(u) Q(g_i)$ as required.  That $Q(g_i)$ is right minimal follows from the same assumption on $g_i$ and the isomorphism $\Hom_{D^b(\rmod \Lambda)}(S_i[-1],Y_i) \cong \stHom_{\Lambda}(\Omega S_i, Y_i)$.  The triangle (5.4) now shows that $Q(\mu^+_{\X}(S_i))[-1] \cong \mu^+_{\X}(Q(S_i))$.   $\Box$ \\

\vspace{2mm}
\noindent
{\it Proof of Okuyama's Lemma.}  By construction and the above results, $$\ul{F}(Q_{\Gamma}(\s')) = Q_{\Lambda}(F(\s')) = Q_{\Lambda}(\mu^+_{\X}(\s)[-1]) = Q_{\Lambda}(\mu^+_{\X}(\s))[-1] = \mu^+_{\X}(Q_{\Lambda}(\s)). \ \  \Box$$

\section{Examples}

We conclude with several remarks regarding simple-minded systems, illustrated by examples, and with an emphasis on a comparison between simple-minded systems and simple-minded collections.  The guiding philosophy behind these two notions is that simple-minded systems resemble the set of simple modules inside the stable category of a self-injective algebra, whereas simple-minded collections resemble the set of simple modules (considered as stalk complexes concentrated in degree 0) inside the derived category of an algebra.  

\vspace{2mm}
\noindent
{\bf Example 1.}  {\it A given triangulated category $\T$ may contain no bricks and hence no simple-minded systems.}  Let $\T$ be any $0$-Calabi-Yau triangulated $k$-category, where $k$ is an algebraically closed field.  The Calabi-Yau condition means that there are trace maps $t_X : \T(X,X) \rightarrow k$ for each $X \in \T$ that yield nondegenerate pairings $\T(X,Y) \times \T(Y,X) \rightarrow k$ via $(f,g) \mapsto t_X(g \circ f) = t_Y(f \circ g)$ for all $X,Y \in \T$ \cite{K3CY}.  Clearly, if $\T(X,X) \cong k$ then any nonzero map $\T(X,Y)$ would be forced to be invertible.  In other words, any brick $X$ must have no nonzero maps to other indecomposables (not isomorphic to $X$) of $\T$.  This condition makes it easy to check that a given $\T$ has no bricks. For example, let $\T= \ul{\mbox{CM}}(R)$ be the stable category of maximal Cohen-Macaulay modules over a simple curve singularity $R$ of Dynkin type over an algebraically closed field $k$ (for instance, take $R = k[[x,y]]/(x^2+y^{n+1})$ of type $\mathbb{A}_n$).  It is well-known that this category is $0$-Calabi-Yau (see 8.3 in \cite{IY}). That no brick $X$ exists in $\T$ can now be seen from the AR-quiver of $\T$, as described in \cite{Yosh}. \\

\vspace{2mm}
\noindent
{\bf Example 2.}  {\it Simple-minded systems in derived categories.}  While stable module categories are the natural setting in which we are interested in simple-minded systems, we may likewise encounter them in derived categories.  For example, if $A$ is a finite-dimensional algebra of finite global dimension, then Happel has proved that there is an equivalence of triangulated categories $D^b(\rmod A) \approx \stgr TA$, where $TA = A \oplus DA$ is the trivial extension algebra, graded by placing $A$ in degree $0$ and $DA$ in degree $1$, and $\stgr TA$ is the stable category of $\mathbb{Z}$-graded $TA$-modules \cite{TCRTA}.  Then one example of a simple-minded system in $\stgr TA$ is given by the (infinite!) set $\s = \{S_i(j)\ |\ 1 \leq i \leq n,\ j\in \mathbb{Z}\}$ with $S_1, \ldots, S_n$ denoting the simple $TA$-modules (concentrated in degree $0$) and $-(j)$ denoting the grading shift.   To understand the corresponding set inside $D^b(\rmod A)$, observe that $S_i \in \stgr TA$ corresponds to the stalk complex consisting of the simple $A$-module $S_i$ concentrated in degree $0$ (we identify the simple modules over $A$ and $TA$).  Furthermore, the grading shift $-(1)$ coincides with the Nakayama functor $\N$ on $\gr TA$, as $\N \cong -\otimes_{TA} D(TA) \cong -\otimes_{TA} TA(1) \cong -(1)$.  Hence, $S_i(j) \in \stgr TA$ corresponds to $\N^j(S_i) \cong \nu^j(S_i[j]) \in D^b(\rmod A)$, where $\nu$ denotes the Serre functor of  $D^b(\rmod A)$.  Moreover, we see that this simple-minded system is invariant under $\nu\circ [1]$ and the subsets at which we can mutate are exactly the $\nu\circ[1]$-orbits of subsets of the simple $A$-modules.  By contrast, observe that any simple-minded {\it collection} in $D^b(\rmod A)$ is necessarily finite, of the same cardinality as the number of simple $A$-modules, up to isomorphism, by Corollary 6.6 in \cite{KoYa}.\\

\vspace{2mm}
\noindent
{\bf Example 3.}  {\it Simple-minded collections in $D^b(\rmod \Lambda)$ do not necessarily induce simple-minded systems in $\stmod \Lambda$.}  As we saw in the last section, if $\Lambda$ is a self-injective algebra, the localization functor $Q : D^b(\rmod \Lambda) \rightarrow \stmod \Lambda$ may take a simple-minded collection to a simple-minded system.  However, this is not always true.  Simple-minded collections, like silting objects, can be mutated with respect to any subset, whereas simple-minded systems may be mutated only with respect to subsets which are stable under $\N$.  Thus to produce an easy counterexample we let $\Lambda$ be the Nakayama algebra whose quiver consists of one 3-cycle and has relations generated by all paths of length 3.  Mutating the simple-minded collection $\{S_1, S_2, S_3\}$ consisting of the simple $\Lambda$-modules at $S_1$ yields $\{S_1[1], S_2, M\}$ where $M$ is isomorphic to the uniserial module of length $2$ with top $S_3$ and socle $S_1$ (as usual, we identify $\rmod \Lambda$ with the stalk complexes concentrated in degree $0$ inside $D^b(\rmod \Lambda)$).  When we pass down to $\stmod \Lambda$ we obtain $\{\Omega^{-1}(S_1), S_2, M\}$ and it is clear that this is not a simple-minded system since $\Hom_{\Lambda}(\Omega^{-1}(S_1), S_2) \neq 0$. \\

It is also interesting to ask whether a given simple-minded system in $\stmod \Lambda$ can be lifted to a simple-minded collection in $D^b(\rmod \Lambda)$ (i.e., is the former the image under the localization functor of some simple-minded collection).  In general, this question has a negative answer, as it is intricately linked to the problem of lifting a stable equivalence to a derived equivalence by work of Linckelmann, Okuyama and Rickard.  We give a brief overview of this connection in order to provide some examples.

    Let us assume $\Lambda$ and $\Gamma$ are indecomposable, non-semisimple, self-injective algebras, and suppose $\alpha : \stmod \Gamma \stackrel{\approx}{\longrightarrow} \stmod \Lambda$ is an equivalence induced by tensoring with a $(\Gamma,\Lambda)$-bimodule (i.e., a stable equivalence of Morita type).  If $\s'$ denotes the set of simple $\Gamma$-modules, then $\alpha(\s')$ is clearly a simple-minded system in $\stmod \Lambda$ which is stable under $\N$.   If $\alpha(\s') = Q(\X)$ for a simple-minded collection $\X \subset D^b(\rmod \Lambda)$, which must also be stable under $\N$, % where $Q : D^b(\rmod \Lambda) \rightarrow \stmod \Lambda$ is the localization functor as in the last section, 
then Al-Nofayee \cite{AlNo} has generalized a theorem of Rickard \cite{EDCSA} to show that there is a derived equivalence $F : D^b(\rmod \Gamma') \rightarrow D^b(\rmod \Lambda)$ sending the simple $\Gamma'$-modules to $\X$.  If we compose the induced stable equivalence with a quasi-inverse of $\alpha$, we obtain a stable equivalence of Morita type from $\Gamma'$ to $\Gamma$ that sends simples to simples.  By a result of Linckelmann, $\Gamma'$ and $\Gamma$ are Morita equivalent \cite{Lin}.  It follows that $\Gamma$ must be derived equivalent to $\Lambda$.  Consequently, whenever we have a stable equivalence $\alpha$ as above between two algebras that are not derived equivalent, the simple-minded system $\alpha(\s') \subset \stmod \Lambda$ cannot be lifted to any simple-minded collection in $D^b(\rmod \Lambda)$.

\vspace{2mm}  \noindent
{\bf Example 4.}  {\it Mutation of simple-minded systems is not necessarily transitive.}  Another interesting consequence of these observations concerns the (lack of) transitivity of mutations of simple-minded systems.  Specifically, in the above setting the images under $\alpha$ of the simple $\Gamma$-modules will provide a simple-minded system in $\stmod \Lambda$ that is not connected to the set of simple $\Lambda$-modules by any sequence of mutations and suspensions.  This follows from Okuyama's Lemma, which implies that any simple-minded system obtained via mutations and suspensions from the set of simple $\Lambda$-modules will coincide with the images of the simple $\Gamma'$-modules under a stable equivalence $\beta: \stmod \Gamma' \rightarrow \stmod \Lambda$ that is induced by a derived equivalence between these two algebras.  As above, if the images of the simple modules under $\alpha$ and $\beta$ coincide in $\stmod \Lambda$, it would follow that $\Gamma$ and $\Gamma'$ are Morita equivalent, and hence that $\Lambda$ and $\Gamma$ are derived equivalent.  

While examples of self-injective algebras which are stably equivalent of Morita type but not derived equivalent are hard to come by, they certainly do exist, even in the context of blocks of group algebras.  For a simpler example that is similar in spirit, we consider stable auto-equivalences given by endo-trivial modules.  Let $k$ be an algebraically closed field of prime characteristic $p$, and let $G$ be a finite $p$-group.  Assume that $L \in \rmod kG$ is an endo-trivial module that is not a (co)syzygy of the trivial module $k$.  Then $kG$ is a local algebra and $-\otimes_k L$ induces an auto-equivalence of $\stmod kG$, sending the unique simple module $k$ to $L$.  It follows that $\{L\}$ is a simple-minded system, which is not a suspension of $\{k\}$, and no mutations are possible here because these sets are singletons.  Furthermore, we also remark that the simple-minded sysem $\{L\}$ cannot be lifted to a simple-minded collection $\{X\}$ in $D^b(\rmod kG)$.  This can be seen as a result of Theorem 5.1 in \cite{EDCSA}, which associates a tilting complex $T$ to such a collection $\{X\}$, should one exist.  But $kG$ is local and symmetric, so the only possibility for $T$ is $kG[i]$ for some $i \in \mathbb{Z}$.   Then property (5.1) from \cite{EDCSA} relating $T$ to $X$ implies that the homology of $X$ is one-dimensional, and it follows that $X$ must be some shift of the trivial module.  But this would imply that $L$ is isomorphic to a (co)syzygy of the trivial module, a contradiction.

 To make this example still more concrete, we can take $G =D_8$ to be the Dihedral group of order $8$.  Then $kG \cong k\langle a, b \rangle/( a^2, b^2, (ab)^2-(ba)^2)$ and $L = kG/ (kG\cdot a + kG \cdot bab)$ is a $3$-dimensional endo-trivial module, which is not a (co)syzygy of the trivial module $k$ (cf. Section 5 in \cite{CarThe}).\\
%Since mutation preserves the cardinality of a simple-minded system, one way of showing that all simple-minded systems in $\stmod \Lambda$ have the same size would be to prove that any two such are connected by a string of mutations and suspensions.   However, simple examples show that this is not true in general.  It would be extremely interesting to identify any cases in which it does hold.

\vspace{2mm} \noindent
{\bf Example 5.} {\it An example of Brou\'e's conjecture.} The following example is taken from Okuyama's preprint \cite{Oku}.  We include it here to illustrate mutation of simple-minded systems.  Let $G$ be the alternating group $A_6$,  $P = \mbox{Syl}_3(G) \cong (\mathbb{Z}/3)^2$ and $H = N_G(P)$.  Let $\Lambda'$ and $ \Lambda $ be the principal blocks of $kG$ and $kH$ respectively, where $k$ is an algebraically closed field of characteristic $3$.   Then $\Lambda$ has the following quiver and its indecomposable projective left modules have the following graphs, with $k$ used as the label for the trivial module and the corresponding vertex of the quiver.
$$\begin{array}{cccccc}
\xymatrixrowsep{1pc} \xymatrixcolsep{1pc}  \xymatrix{k \ar[r]<.5ex> \ar[d]<.5ex> & 1 \ar[d]<.5ex> \ar[l]<.5ex> \\ 3 \ar[u]<.5ex> \ar[r]<.5ex> & 2 \ar[u]<.5ex> \ar[l]<.5ex>} & \hspace{3mm} &  
\small{\xymatrixrowsep{-.3pc} \xymatrixcolsep{-.3pc} \xymatrix{ & & k \\ & 1 & & 3 \\ 2 & & k & & 2 \\ & 1 & & 3 \\ & & k}} & 
\small{\xymatrixrowsep{-.3pc} \xymatrixcolsep{-.3pc} \xymatrix{& & 1 \\ & 2 & & k \\ 3 & & 1 & & 3 \\ & 2 & & k \\ & & 1}} &
\small{\xymatrixrowsep{-.3pc} \xymatrixcolsep{-.3pc} \xymatrix{& & 2 \\ & 3 & & 1 \\ k & & 2 & & k \\ & 3 & & 1 \\ & & 2}} &
\small{\xymatrixrowsep{-.3pc} \xymatrixcolsep{-.3pc} \xymatrix{& & 3 \\ & k & & 2 \\ 1 & & 3 & & 1 \\ & k & & 2 \\ & & 3}}
\end{array}$$

It is known that restriction from $G$ to $H$ induces an equivalence $\stmod \Lambda' \rightarrow \stmod \Lambda$.  We write $Z_i\  (0 \leq i \leq 3)$ for the restrictions of the simple $\Lambda'$-modules, which have graphs as follows.

$$\begin{array}{cccc} 
\small{Z_0 = k} & 
\small{Z_1 =  \vcenter{ \xymatrixrowsep{-.3pc} \xymatrix{1 \\ k \\ 3}} }&
\small{Z_2 = \vcenter{\xymatrixrowsep{-.3pc} \xymatrixcolsep{-.3pc} \xymatrix{& 2 \\  1 & & 3 \\  & 2 }} }&
\small{Z_3 =  \vcenter{\xymatrixrowsep{-.3pc} \xymatrix{3 \\ k \\ 1}}}
\end{array}$$

By performing two mutations we can transform this simple-minded system into the set of simple $\Lambda$-modules.  First, we apply $\mu^+_{\{Z_0, Z_2\}}$.  This leaves $Z_0$ and $Z_2$ unchanged and takes $Z_1$ to  $S_1$ and $Z_3$ to $S_3$ by virtue of the triangles $S_1 \rightarrow \Omega Z_1 \rightarrow X \rightarrow$  and $S_3 \rightarrow \Omega Z_3 \rightarrow Y$, where $X$ (resp. $Y$) is an extension of $Z_0$ by $Z_2$ and hence in $\F(\{Z_0,Z_2\})$.  Next, we perform $\mu^+_{\{k,S_1,S_3\}}$, which takes $Z_2$ to $S_2$ by virtue of the triangle $S_2 \rightarrow \Omega Z_2 \rightarrow W \rightarrow$ where $W \in \F(\{k,S_1,S_3\})$ is the direct sum of two uniserial modules with top $k$ and with socles isomorphic to $S_1$ and $S_3$ respectively.  Therefore, according to the above discussion, there is a derived equivalence $D^b(\rmod \Lambda') \stackrel{\approx}{\rightarrow} D^b(\rmod \Lambda)$.

\end{document}